\documentclass[12pt]{article}
\usepackage{amsmath}
\usepackage{amssymb}
\usepackage{amsthm}
\usepackage{amscd}
\usepackage{amsfonts}
\usepackage{graphicx}%
\usepackage{fancyhdr}

\theoremstyle{plain} \numberwithin{equation}{section}
\newtheorem{theorem}{Theorem}[section]
\newtheorem{corollary}[theorem]{Corollary}

\newtheorem{lemma}[theorem]{Lemma}

\theoremstyle{definition}

\newcommand{\A}{\mathcal{A}}
\newcommand{\B}{\mathcal{B}}
\newcommand{\C}{\mathcal{C}}
\newcommand{\R}{\mathcal{R}}
\newcommand{\h}{\mathcal{H}}

\begin{document}
\title{A description of the logmodular subalgebras in the finite dimensional $C^*$-algebras.}

\author{Kate Juschenko}
\date{}

\maketitle

\begin{abstract}We show that every logmodular subalgebra of $M_n(\mathbb{C})$ is unitary equivalent to an algebra of block upper triangular matrices, which was conjectured in \cite{VM}. In particular, this shows that every unital contractive representation of a logmodular subalgebra of $M_n(\mathbb{C})$ is automatically completely contractive.
\end{abstract}
\section{Introduction}

Let $\B$ be a unital $C^*$-algebra and $\A\subseteq \B$ be a closed unital subalgebra. Denote by $\A^{-1}$ the set of invertible  elements in $\A$. A non-commutative generalization of logmodularity of closed subalgebras of continuous functions on a compact Hausdorff space was introduced  in \cite{BL}, see also \cite{BL2}. Namely, $\A$ is \textit{logmodular} in $\B$ if the set $\{a^*a:a\in\A^{-1}\}$ is dense in the set of positive invertible elements of $\B$.

In \cite{VM}, V. Paulsen and M. Raghupathi study the question when a contractive homomorphism of a logmodular algebras into $B(\h)$ is completely contractive, where $B(\h)$ denotes the algebra of all bounded operators on a Hilbert space $\h$. Generalizing C. Foias and I. Suciu, \cite{FS}, they showed that a contractive homomorphism $\pi:\A\rightarrow B(\h)$ is completely contractive if and only if its second amplification $\pi^{(2)}=\pi\otimes 1_2:\A\otimes M_2(\mathbb{C})\rightarrow B(\h)\otimes M_2(\mathbb{C})$ is contractive. It is still unknown if there are logmodular algebras with contractive but not completely contractive homomorphisms. In particular, it is interesting to decide this question for $H^{\infty}(\mathbb{D})$ considered as logmodular subalgebra of $L^{\infty}(\mathbb{T})$ and the algebra of upper triangular matrices on infinite dimensional Hilbert space.

The decomposition of Cholesky shows that an algebra of block upper triangular matrices is logmodular in $M_n(\mathbb{C})$. It was proved by V. Paulsen and M. Raghupathi in \cite{VM} that if $\A\subseteq M_n(\mathbb{C})$ is logmodular and contains the diagonal matrices $D_n$ then it is unitary equivalent to an algebra block upper triangular matrices. It was conjectured that algebras of block upper triangular matrices are the only subalgebras on $M_n(\mathbb{C})$ that have the logmodularity property. In Theorem \ref{main} we prove this in the affirmative. It was proved in \cite{PP} that all contractive homomorphisms of an algebra of block upper triangular matrices are completely contractive. Thus there are no logmodular subalgebras in $M_n(\mathbb{C})$ that can provide us with examples of contractive but not completely contractive homomorphisms. 

 We will use the following notations. Matrix units of $M_n(\mathbb{C})$ will be denoted by $E_{i,j}$. Put $P_{\{i_1,..,i_n\}}= \sum \limits_{i\in\{i_1,..,i_n\}} E_{i,i}$. Given matrix $a=[a_{ij}]\in M_n(\mathbb{C})$ we denote $a(i,j)=a_{ij}$.

\section{A description of the logmodular subalgebras in $M_n(\mathbb{C})$.}

The proof of the main result will be divided into several lemmas.
Note that if $\A$ is logmodular in $\B$ then $\A^*$ is also logmodular in $\B$. Assume $\A\subseteq M_n(\mathbb{C})$ is logmodular, then $\A$ is unital and by compactness argument for every positive $b\in M_n(\mathbb{C})$ there exist $a,c\in\A$ such that $b=a^*a=cc^*$. In particular we have the following lemma.

\begin{lemma}\label{l1}
If $\mathcal{A}$ is logmodular in $M_n(\mathbb{C})$ then there are $\alpha_{ij}$, $\beta_{ij}\in \mathbb{C}$, $i,j\in\{1,..,n\}$ such that each row of $[\alpha_{ij}]$ and each column of $[\beta_{ij}]$ is non-zero and
$\sum\limits_{k} \alpha_{ik}E_{i,k}\in\mathcal{A}$, $\sum\limits_{k} \beta_{kj}E_{k,j}\in\mathcal{A}$ for every $i,j\in\{1,..,n\}$
\end{lemma}
\begin{proof}
By logmodularity of $\A$ we have $E_{i,i}=a_ia_i^*=b_i^*b_i$ for some $a_i, b_i\in\A$. Then we can put $\alpha_{ij}=a_i(i,j)$ and $\beta_{ij}=b_i(i,j)$ satisfy the statement.
\end{proof}

\begin{lemma}\label{l2}
Let $v=(v_1,..,v_{n})\in \mathbb{C}^n$, $||v||=1$ and
\begin{gather*}
A=\text{span }(\sum\limits_{j=1}^n v_j E_{1,j},..,\sum\limits_{j}v_j E_{n,j})\subset M_n(\mathbb{C})\\
B=\text{span }(\sum\limits_{i=1}^n v_i E_{i,1},..,\sum\limits_{i}v_i E_{i,n})\subset M_n(\mathbb{C})
\end{gather*}
Then there exist unitary $U_1, U_2\in M_{n}(\mathbb{C})$ such that $E_{1,1}\in U_1\A U_1^*$ and $E_{1,1}\in U_2 \B U_2^*$.
\end{lemma}
\begin{proof}
Let $U_1\in M_n(\mathbb{C})$ be unitary such that $vU_1=(1,0,..,0)$ then
\begin{gather*}
U_1^*A U_1=\{U_1^*\left(
       \begin{array}{ccc}
         \alpha_1v_1 & .. & \alpha_1v_n \\
         .. & .. & .. \\
         \alpha_nv_1 & .. & \alpha_nv_n \\
       \end{array}
     \right)U_1: \alpha_1,\ldots, \alpha_n\in\mathbb{C}\}= \\=\{U_1^*\left(
                      \begin{array}{c}
                        \alpha_1vU \\
                        .. \\
                        \alpha_nvU \\
                      \end{array}
                    \right):\alpha_1,\ldots, \alpha_n\in\mathbb{C}\}=\\=\{ U_1^*\left(
                                    \begin{array}{cccc}
                                      \alpha_1 & 0 & .. & 0 \\
                                      .. & .. &  & .. \\
                                      \alpha_n & 0 & .. & 0 \\
                                    \end{array}
                                  \right):\alpha_1,\ldots, \alpha_n\in\mathbb{C}\}\subseteq M_n(\mathbb{C})E_{1,1}.
\end{gather*}
Since $\dim(A)=n$ we have $U_1^*A U_1=M_n(\mathbb{C})E_{1,1}$ and therefore $E_{1,1}\in U_1^*A U_1$. The same argument applied to $B^*$ provides the existence of $U_2$.
\end{proof}

In order to state the next lemma we need some definitions from general graph theory. A directed graph $(V,\R)$ is {\it transitive} if $(i,j)\in \R$ and $(j,k)\in \R$ implies $(i,k)\in R$ for every $i,j,k\in V$. A directed graph $(V,\R)$ is {\it full graph} if $(i,j)\in \R$ for every $i,j\in V$. We will say that a subgraph $(V_1,\R_1)$ of $(V,\R)$ {\it does not have an exit} if conditions $i\in V_1$, $j\in V$ and $(i,j)\in \R$ imply $j\in V_1$. A vertex $i\in V$ has {\it an exit} if there is $j\in V, j\neq i$ such that $(i,j)\in \R$. 

\begin{lemma}\label{graph}
Let $(V,\R)$ be a directed and finite graph. Assume that the graph $(V,\R)$ is transitive and every $i\in V$ has an exit. Then $(V,\R)$ contains a full subgraph $(V_1,\R_1)$ which does not have an exit.
\end{lemma}
\begin{proof}
Note that by transitivity of $\R$ for every cycle $\C$ in $(V,\R)$ the subgraph induced on $\C$ by $\R$ is automatically a full graph.
Let $i_0\in V$ be an arbitrary vertex. Then there is an exit $i_1\in V$ with $(i_0,i_1)\in \R$. Let $i_2\in V$ be an exit of $i_1$. Continuing this process we obtain a directed path with vertices $i_0,\ldots, i_n\in V$. Since the graph is finite we will have a cycle in our sequence. Assume that this cycle is the maximal in the sense that it is not contained in any other full subgraph of $(V,\R)$. Then if this cycle has a vertex $i$ with an exit, we will repeat the procedure with $i_0=i$. Note that we can not return to the previous sequences since it will contradict the maximality of the cycles that have been obtained before. Since the graph is finite we will eventually reach a cycle without any exit.
\end{proof}

\begin{lemma}\label{crutial}
Let $v_1,..,v_n$ be non-zero vectors in $\mathbb{C}^n$ and $\A$ be the algebra generated by
$$A_1=\sum\limits_{j=1}^n v_1(j) E_{1,j},\ldots, A_n=\sum\limits_{j=1}^n v_n(j) E_{n,j}$$
Then there exist $k\in\{1,..,n\}$, a non-zero $v\in\mathbb{C}^k$ and unitary $U\in M_n(\mathbb{C})$
such that
$$\sum\limits_{j=1}^k v(j) E_{1,j},\ldots, \sum \limits_{j=1}^k v(j) E_{k,j}\in U\A U^*$$
\end{lemma}
\begin{proof}
Denote by $\R\subseteq\{1,..,n\}\times\{1,..,n\}$ the set of indices with the property $(i,j)\in \R$ iff
$S(i,j)\neq 0$ for some $S$ which is a word in the generators $A_1,\ldots, A_n$. Thus if $(i,j)\in \R$
and $(j,k)\in \R$ then $(i,k)\in \R$. Since $v_1,\ldots, v_n$ are all non-zero we have that for every $m\in\{1,..,n\}$ there exists $i\in\{1,..,n\}$ such that $(m,i)\in \R$. 

Consider the graph $(\{1,..,n\},\R)$.  Let $\R_1$ be the set of edges of a connected component of the directed graph $(\{1,..,n\}, \R)$ and let $V$ be the set of vertices of $\R_1$. Then Lemma \ref{graph} implies that $\R_1$ contains a full graph without exit, denote the set of its vertices by $W=\{i_1,..,i_k\}$ and the set of its edges by $\R'$. Chose a unitary $U$ that maps $e_{i_1},\ldots,e_{i_k}$ to $e_1,\ldots,e_k$ by permutation of the basis elements.

Since $\R'$ does not have any exit we have that for every $(i,j)\in \R'$ there exists a word $S$ in the generators $A_1,\ldots, A_n$ such that $S(i,j)\neq 0$ and $S\in P_{\{i_1,..,i_k\}} M_n(\mathbb{C}) P_{\{i_1,..,i_k\}}$.

Since $(W, \R')$ is full for every $t\in\{i_1,..,i_k\}$ and $i\in \{i_1,..,i_k\}$ we have $(i,t)\in\R_1$. Let $S_i$ be the word on the generators $A_1,\ldots, A_n$ that provides $(i,t)\in \R_1$, then $S_i=E_{i,i}S_i$ and $S_i S_t= S_i(i,t) E_{i,t}S_t$. Thus we have the statement of the lemma with $v=(S_t(t,i_1),\ldots, S_t(t,i_k))$.
\end{proof}

\begin{lemma}\label{last}
Fix $k\in\{1,\ldots n\}$ and let $A\subseteq P_{\{1,\ldots, k\}} M_{n}(\mathbb{C})$ be a set
such that for every $m\in\{1,\ldots,k\}$ and any set of indices $\{i_1,\ldots, i_m\}\subseteq \{1,\ldots, k\}$ we have $E_{m,m}\in A$ and $P_{\{i_1,\ldots,i_m\}}A (1_n-P_{\{i_1,\ldots,i_m\}})\neq \{0\}$. Then for every $t\in \{1,\ldots, k\}$ the algebra generated by $A$ contains an  element $S$ such that $S(t,l)\neq 0$ for some $l\in \{k+1,\ldots,n\}$.
\end{lemma}
\begin{proof}
Let $I=\{i:E_{i,i}\B(1_n-P_{\{1,\ldots, k\}})=0, 1\leq i\leq k\}$ where $\B$ is the algebra generated by $A$. To reach a contradiction assume that $I\neq \emptyset$.
Permuting the part of the basis $e_1,\ldots, e_k$ we can assume $I=\{1,\ldots, d\}$. Since $P_{\{1,\ldots, k\}}A (1-P_{\{1,\ldots, k\}})\neq 0$ we have $d<k$.  Take $T\in P_{\{1,\ldots,d\}} A (1_n-P_{\{1,\ldots,d\}})$ and $T(i,j)\neq 0$ for some $1\leq i\leq d$, $d< j \leq k$. There exists $P\in E_{j,j}\B (1-P_{\{1,\ldots,k\}})$ with $P(j,l)\neq 0$ for some $k< l\leq n$. Then the $(i,l)$-th entry of $E_{i,i}TE_{j,j}P\in E_{i,i}\B (1-P_{1,\ldots,k})$ is non-zero which contradicts $i\in I$.
\end{proof}

Now we are ready to prove our main result.

\begin{theorem}\label{main} If $\A\subseteq M_n(\mathbb{C})$ is logmodular then $\A$ is unitary equivalent to  an algebra block upper triangular matrices.
\end{theorem}
\begin{proof}
We will proceed by induction on dimension $n$. For $n=1$ the statement is trivial. Assume that all logmodular subalgebras in $M_k(\mathbb{C})$, $k<n$, are unitary equivalent to block upper triangular matrices. Let $\A$ be logmodular in $M_n(\mathbb{C})$. 

By Lemma \ref{l1} we have that there are non-zero $v_1,\ldots, v_n\in \mathbb{C}^n$ such that $\sum \limits_{j=1}^{n}v_1(j) E_{1,j},\ldots \sum \limits_{j=1}^{n}v_n(j) E_{n,j} \in \A$. Then by Lemma \ref{crutial} there exist $k\in \{1,\ldots, n\}$, $v\in \mathbb{C}^k$ with $||v||=1$ and a unitary $V\in M_n(\mathbb{C})$ such that $\sum \limits_{j=1}^{k}v(j) E_{1,j},\ldots, \sum \limits_{j=1}^{k}v(j) E_{k,j}\in V\A V^*$. Thus by Lemma \ref{l2} we have $E_{1,1}\in V\A V^*$ for some unitary $U\in M_n(\mathbb{C})$. \\

We will prove by induction that $E_{1,1},\ldots, E_{n,n}\in U\A U^*$ for some unitary $U \in M_n(\mathbb{C})$.
Assume that $E_{1,1},\ldots, E_{k,k}\subseteq V\A V^*$ for some $V\in \mathcal{U}_n(\mathbb{C})$ and $k<n$. Denote $V\A V^*$ again by $\A$.

Firstly, assume the existence of the set $\{i_1,\ldots,i_m\}\subseteq \{1,\ldots, k\}$ such that $P_{\{i_1,\ldots, i_m\}} a (1_n-P_{\{i_1,\ldots, i_m\}})=0$ for every $a\in \A$. Denote $\B=(1_n-P_{\{i_1,\ldots, i_m\}})\A(1_n-P_{\{i_1,\ldots, i_m\}})$ and $\C=(1_n-P_{\{i_1,\ldots, i_m\}})M_n(\mathbb{C})(1_n-P_{\{i_1,\ldots, i_m\}})$. Then $\B$ is logmodular in $\C$. Indeed, let $a\in \C_{+}$ then there exists $b\in \A$ such that $a=b^*b$, but $P=(1_n-P_{\{i_1,\ldots, i_m\}})\in \A$, therefore $a=(bP)^*bP$ and $bP\in \B$. Since $\C$ is unitary equivalent to $M_t(\mathbb{C})$ for some $t<n$ we have that $\B$ is unitary equivalent to block upper triangular matrices and thus $D_n\subseteq U\A U^*$ for some unitary $U$.

Thus we arrived to the case that for every $m\in\{1,..,k\}$ and every subset $\{i_1,\ldots, i_m\}\subseteq \{1,\ldots, k\}$ there exist an element $a\in \A$ such that $P_{\{i_1,\ldots, i_m\}}a(1_n-P_{\{i_1,\ldots, i_m\}})\neq 0$.
We claim that for $P=P_{\{k+1,\ldots,n\}}=1_n-P_{\{1,\ldots, k\}}\in \A$ the set $P\A P\subset \A$ contains
\begin{gather*}
A_1=\sum\limits_{j=k+1}^{n} v_1(j)E_{k+1,j},\ldots,
A_{n-k}=\sum\limits_{j=k+1}^{n} v_{n-k}(j)E_{n,j}
\end{gather*}
for some non-zero $v_i\in \mathbb{C}^{n-k}$, $1\leq i\leq n-k$. Then from the claim, Lemma \ref{l2} and Lemma \ref{crutial} follows that there exists a unitary $U\in M_{n-k}(\mathbb{C})$ such that $$E_{k+1,k+1}\in (1_k\oplus U)P\A P (1_k\oplus U^*).$$
Therefore $(1_k\oplus U)\A  (1_k\oplus U^*)$ contains $E_{1,1},\ldots, E_{k+1,k+1}$ and by induction we have the statement of the theorem.\\

To prove the claim consider a decomposition of $E_{t,t}$ in $\A$ for $k+1\leq t\leq n$, that is $R_t R_t^*=E_{t,t}$, where $R_t\in\A$ is a matrix with all rows equal to zero except for the $t$-th row. We will find $v_i\in \mathbb{C}^{n-k}$ by "shifting" non-zero elements from the set $R_t P_{\{1,\ldots,k\}}$ to the set $P\A P$.
Note that, now $P_{\{1,\ldots,k\}}\A\subseteq \A$ satisfies the assumptions of Lemma \ref{last}.

 Assume that $R_t(t,i)\neq 0$ for some $1\leq i\leq k$ and $k+1\leq t\leq n$. By Lemma \ref{last} there exists $S\in \A$ with $S(i,t)\neq 0$. Denote $R_t E_{i,i} S$ again by $R_t$. Doing this process with all $R_t$, $k+1\leq t\leq n$ we obtain a set of non-zero rows with the property $(1_n-P_{\{1,\ldots,k\}})R_t (1_n-P_{\{1,\ldots,k\}})=R_t\neq 0$ for every $k+1\leq t\leq n$. Then the vectors $v_i$ with $v_i(j)=R_{i+k}(i+k,j)$, $1\leq i\leq n-k$ have the required property.

\end{proof}
As a consequence of the Theorem \ref{main} and the fact that all contractive homomorphisms of an algebra of block upper triangular matrixes are completely contractive, see \cite{PP}, we get the following corollary.
\begin{corollary}
If $\A$ is a logmodular subalgebra in $M_n(\mathbb{C})$ then every contractive unital homomorphism $\pi:\A\rightarrow B(H)$ is completely contractive.
\end{corollary}

\noindent{\bf Acknowledgement.}
The author is grateful to Gilles Pisier for suggesting a number of corrections and for constant encouraging.
The work was partially supported by NSF grant 0503688.

\end{document}